\documentclass[12pt,reqno]{amsart}

\usepackage{latexsym}
\usepackage[centertags]{amsmath}
\usepackage{amsfonts}
\usepackage{amssymb}
\usepackage{amsthm}
\usepackage{newlfont}
\usepackage[square,sort]{natbib}

\def\mytopsep{3mm}

\newtheoremstyle{myplain}{\mytopsep}{\mytopsep}{\itshape}{0pt}{\bfseries}{.}{3mm}{}
\newtheoremstyle{mydefinition}{\mytopsep}{\mytopsep}{\normalfont}{0pt}{\bfseries}{.}{3mm}{}

\newtheoremstyle{myremark}{\mytopsep}{\mytopsep}{\normalfont}{0pt}{\bfseries}{.}{3mm}{}

\theoremstyle{myplain}

\newtheorem{thm}{Theorem}[section]
\newtheorem{cor}[thm]{Corollary}
\newtheorem{lem}[thm]{Lemma}

\theoremstyle{mydefinition}

\theoremstyle{myremark}
\newtheorem{rem}[thm]{Remark}

\allowdisplaybreaks[1]

\makeatletter\@addtoreset{equation}{section}\makeatother

\def\diag{\mathrm{diag }}

\def\NN{\mathbb{N}}

\def\ct{\mathop{\mathrm{CT}}}

\begin{document}

\title{Constructing All Magic Squares of Order Three}

\author{Guoce Xin}
\address{Department
of Mathematics\\
Brandeis University\\
Waltham MA 02454-9110} \email{guoce.xin@gmail.com}

\date{September 23, 2004}
\begin{abstract}

We find by applying MacMahon's partition analysis that all magic
squares of order three, up to rotations and reflections, are of
two types, each generated by three basis elements. A combinatorial
proof of this fact is given.
\end{abstract} \maketitle

{\small \bf Keywords:} {\small \it magic square, linear
Diophantine equations}

\section{introduction}
A \emph{magic square} of \emph{order} $n$ is an $n$ by $n$ matrix
with distinct nonnegative integer entries such that every row sum,
column sum, and (two) diagonal sums equals to the same number $m$,
the \emph{magic number}. Adding $1$ to every entry will give us a
traditional magic square of positive integers. A magic square is
\emph{pure} if the entries are the consecutive numbers from $0$ to
$n^2-1$, and hence it has magic number $3\binom{n+1}{3}$.
%A \emph{semimagic square} satisfies all the restrictions of a magic
%square except the conditions on the diagonals.

Magic squares have been objects of study for centuries. As
Pickover wrote in his book% ``The Zen of Magic Squares, Circles,
%and Stars"
 \cite[p. 60]{pickover}:

\vspace{4mm}
\parbox{16cm}{\emph{
\dots the holy grail of magic squares creation would be to
discover a method that would generate every possible arrangement
for a square of a given size. Such a solution is probably not
discoverable}.}

\vspace{4mm} This ``holy grail" could be achieved by first finding
the complete generating function (which is a rational function)
for magic squares of a given size, and then writing the generating
function as a sum of simple rational functions, the series
expansion of which has only nonnegative coefficients.

We achieve this for magic squares of order $3$, as given in
Theorem \ref{t-main}.

\emph{Weak magic squares}, magic squares without the restriction
of distinct elements, have been studied in
\cite{beck_magicsquare,andrews5,ahmed_magicsquare,ahmed_magiccubes}
by using the rich theory of counting solutions of a system of
linear Diophantine equations, or equivalently, counting lattice
points of a convex polytope. For further references, see
\cite[Ch.~4.6]{EC1}. These methods also apply to counting magic
squares, but give no obvious reason why a simple solution as in
Theorem \ref{t-main} exists.

We give our main result in Section 2, and give a combinatorial
proof in Section 3. In Section 4, we discuss the discovery of our
main result and possible future work.

\section{Main Results}
A magic square of order $3$ is a $3$ by $3$ matrix of distinct
nonnegative integers such that every row sum, column sum, and
diagonal sum equals the magic number $m$.

Our main result is the following Theorem \ref{t-main}, which
generates all magic square of order $3$.

Let
\begin{equation}\label{e-abcd}
A=\left[ \begin {array} {ccc} 1&1&1\\1&1&1\\1&1&1
\end {array} \right],\;\;
B=\left[ \begin {array}{ccc} 5&0&4\\2&3&4
\\2&6&1\end {array} \right], \;\;
C=\left[ \begin {array}{ccc} 2&0&1\\0&1&2
\\1&2&0\end {array} \right], \;\;
D=\left[ \begin {array}{ccc} 3&0&3
\\2&2&2\\1&4&1\end {array}
 \right],\end{equation}
 \begin{equation}
T_1= \left[ \begin {array}{ccc} 7&0&5\\2&4&6
\\3&8&1\end {array} \right],\qquad
T_2=
 \left[ \begin {array}{ccc} 8&0&7\\4&5&6
\\3&10&2\end {array} \right]
.\label{e-type}
\end{equation}

Then they are related as follows:
\begin{equation}
\label{e-depdence} B=C+D, \qquad   T_1=B+C,\qquad T_2=B+D.
\end{equation}
If we let $C'$ be obtained from $C$ by reflecting in the vertical
axis, then we have one more relation: $D=C+C'$. It is
straightforward to check that $A,C,$ and $D$ are linearly
independent.

In fact, $A,C,D$ are the three basis elements that generate all
magic squares of order $3$, and $T_1, T_2$ are the unique magic
squares with magic numbers $12$ and $15$, respectively, up to
rotations and reflections.

\begin{thm}\label{t-main}
Every magic square of order three, up to rotation and reflection,
can be written uniquely as either $T_1+iA+jB+kC$ or $T_2
+iA+jB+kD$, where $i,j,k$ are nonnegative integers and $A,B,C,D,
T_1, T_2$ are as in \eqref{e-abcd}, \eqref{e-type}.
\end{thm}

\begin{rem}
Note that traditional magic squares can be generated by either
$iA+jB+kC$ or $iA+jB+kD$ for {\rm positive} integers $i,j,k$. This
description reveals a kind of symmetry.
\end{rem}

Theorem \ref{t-main} says that magic squares, as a set of lattice
points, is a disjoint union of $16=8\cdot 2$ \emph{polyhedrons}
that are isomorphic to $\NN^3$, where the factor $8$ is the order
of the dihedral group of rotations and reflections. We will give a
combinatorial proof of this result in the next section.
%Polyhedrons are defined by some linear equations.

\begin{cor}
The number of magic squares of order $3$ with magic number $3s$
and its associated generating function is given by
\begin{align*}
\frac{8t^4(1+2t)}{(1-t)(1-t^2)(1-t^3)}&=\sum_{s\ge 0} %\left(4\binom{m+2}{2}-\frac{38}{3}\binom{m+1}{1}+\frac{37}{3}-(-1)^m-(m\bmod 3)-1 \right)
\left(2s^2-\frac{20}{3}s+1-(-1)^s+\frac{8}{3}(s\bmod
3)\right)t^s\\
&=8
\left({t}^{4}+3{t}^{5}+4{t}^{6}+7{t}^{7}+10{t}^{8}+13{t}^{9}+17{
t}^{10}+22{t}^{11}+26{t}^{12}+\cdots \right).
\end{align*}
\end{cor}

\section{A Combinatorial Proof}

In what follows, magic squares are always of order $3$ unless
specified otherwise.

Let $M$ be a magic square with magic number $m$. We write
\begin{align}
M=\left[ \begin {array} {ccc}
a_1&a_2&a_3\\b_1&b_2&b_3\\c_1&c_2&c_3
\end {array} \right],
\end{align}
where
\begin{enumerate}
\item[$\mathbf{C1}$:] Every row sum, column sum, and diagonal sum is equal to $m$.
\item[$\mathbf{C2}$:] The entries of $M$ are distinct nonnegative integers.
\end{enumerate}

Rotating or reflecting $M$ will give us different magic squares.
Without loss of generality, we can assume that $c_3$ is smaller
than $a_1,a_3$ and $c_1$, and that $c_1<a_3$. Also by subtracting
 $A$ times the minimal entry of $M$ from $M$, we can assume that
$0$ is an entry of $M$. Then $M$ satisfies the following two extra
conditions:
\begin{enumerate}
\item[$\mathbf{C3}$:] One of the entries of $M$ is $0$.
\item[$\mathbf{C4}$:] $c_3<a_1,a_3,c_1$, and $c_1<a_3$.
\end{enumerate}

In fact, $\mathbf{C4}$ can be replaced with
\begin{enumerate}
\item[$\mathbf{C4'}$:] $c_3<c_1<a_3<a_1$,
\end{enumerate}
which follows from the sum of the two diagonals.

If $M$ satisfies the above four conditions, then we say that $M$
is a \emph{reduced} magic square. It is well-known that the magic
number $m$ is $m=3b_2$. Let $m=3s$ or equivalently $s=b_2$.
\begin{lem}\label{l-main}
If $M$ is a reduced magic square, then $a_2=0$, and $c_2=2s$.
\end{lem}
\begin{proof}
Since $b_2=s=m/3$, where $m$ is the magic number, the last
statement $c_2=2s$ follows from $a_2=0$, which is what we are
going to show now.

In a reduced magic square $M$, $a_1$ and $a_3$ are the largest two
entries among the four corners $a_1,a_3,c_1,c_3$. It follows from
the first and third row sums and column sums that
$a_2<b_1,c_2,b_3$. We see that $b_2=s\ge 4$, since all entries are
distinct nonnegative integers. It remains to show that $c_3$
cannot be $0$. Assuming that $c_3$ equals $0$, then $a_1=2s$ by
the diagonal sum $a_1+b_2+c_3=3s$. By investigating the first row
sum and the first column sum, we get $a_3<s-1$ and $c_1<s-1$,
contradicting the condition for the diagonal $(a_3, b_2, c_1)$.
\end{proof}

\begin{lem}\label{l-main1}
A reduced magic square $M$ can be uniquely written as $T_1+\alpha
C+\beta D$, where $\alpha\ge -1$ and $\beta\ge 0$ are integers.
\end{lem}
\begin{proof}
To see the existence, we use Lemma \ref{l-main}. Assuming that
$c_3=r$ and $b_2=s$, we obtain all the entries of $M$ by the
condition $\mathbf{C1}$ for row sums, column sums, and diagonal
sums:
$$M=\begin{bmatrix}
  2s-r & 0 & s+r \\
  2r & s & 2s-2r \\
  s-r & 2s & r \\
\end{bmatrix}.$$
Comparing the above matrix with
$$T1+\alpha C+\beta D=\begin{bmatrix}
  7+2\alpha+3\beta & 0 & 5+\alpha+3\beta \\
  2+2\beta & 4+\alpha+2\beta & 6+2\alpha+2\beta \\
  3+\alpha+\beta & 8 +2\alpha+4\beta & 1+\beta \\
\end{bmatrix}, $$we solve uniquely for $\alpha$ and $\beta$:
$$\alpha=s-2r-2, \text{ and } \beta =r-1.$$
Consequently,
$$s= \alpha +2\beta+4 \text{ and } r=\beta +1.$$

We see that $c_3\ge 1$ and $c_1>c_3$ implies that $\alpha\ge -1$
and $\beta\ge 0$, completing the proof of the existence.

The uniqueness follows from the above proof, and also from the
fact that $C$ and $D$ are linearly independent.
\end{proof}
We are now ready to give the proof of our main theorem.

\begin{proof}[Proof of Theorem \ref{t-main}]
It is straightforward to check that \begin{align}\label{e-MT1}
T_1+iA+jB+kC&= \left[
\begin {array}{ccc} 7+i+5\,j+2\,k&i&5+i+4\,j+k
\\2+i+2\,j&4+i+3\,j+k&6+i+4\,j+2\,k
\\3+i+2\,j+k&8+i+6\,j+2\,k&1+i+j\end {array}
 \right],\\
 T_2 +iA+jB+kD &=\left[ \begin {array}{ccc} 8+i+5\,j+3\,k&i&7+i+4\,j+3\,k
\\4+i+2\,j+2\,k&5+i+3\,j+2\,k&6+i+4\,j+2\,k
\\3+i+2\,j+k&10+i+6\,j+4\,k&2+i+j+k\end {array}
 \right]\label{e-MT2}
 \end{align}
give different magic squares for all nonnegative integers $i,j,k$.

Given a magic square $M$, we need to show that $M$ equals either
\eqref{e-MT1} or \eqref{e-MT2}.

Let $i$ be the minimum of the entries of $M$. Then up to rotations
and reflections, we can assume $M'=M-iA$ is a reduced magic
square. By Lemma \ref{l-main1}, $M'$ can be uniquely written as
$T_1+\alpha C+\beta D$, with $\alpha\ge -1$ and $\beta\ge 0$.

If $\alpha\ge \beta \ge 0$, $M'$ can be rewritten (recall that
$B=C+D$) as $T_1+\beta B+(\alpha-\beta)C$. Hence we let
$j=\beta\ge 0$ and $k=\alpha-\beta\ge 0$.

If $\alpha<\beta$, $M'$ can be rewritten (recall that
$T_1+D=T_2+C$) as
$$T_1+\alpha B+(\beta-\alpha)D=T_2+C+\alpha B+(\beta-\alpha-1)D
=T_2+(\alpha+1)B+(\beta-\alpha-2)D.$$ Thus we let $j=\alpha+1\ge
0$ and $k=\beta-\alpha -2\ge -1$.

The only remaining case is $k=-1$, which is equivalent to
$\beta=\alpha+1$. But in this case
$$M'=T_1+(\beta-1)C+\beta D=
\left[ \begin {array}{ccc} 5+5\,\beta&0&4+4\,\beta
\\2+2\,\beta&3+3\,\beta&4+4\,\beta
\\2+2\,\beta&6+6\,\beta&1+\beta\end {array} \right]
,
$$
which is not a magic square because it has equal entries.
\end{proof}

\section{Further Discussion}

The combinatorial proof in the previous section seems unlikely to
be applicable to magic squares of higher order. We describe how we
discovered Theorem \ref{t-main} by using MacMahon's partition
analysis, which has been restudied by Andrews and his coauthors in
a series of papers (see e.g., \cite{andrews5}).

MacMahon's idea is to use new variables to replace linear
constraints. For example, if we want to count nonnegative integral
solutions of the linear equation $a_1+a_2-a_3=0$, we can simply
write the generating function as
$$\sum_{a_1,a_2,a_3\ge 0\atop a_1+a_2-a_3=0}x_1^{a_1} x_2^{a_2} x_3^{a_3}=
\sum_{a_1,a_2,a_3\ge 0} \ct_\lambda \lambda^{a_1+a_2-a_3}x_1^{a_1}
x_2^{a_2} x_3^{a_3}=\ct_\lambda \frac{1}{(1-\lambda x_1)(1-\lambda
x_2 ) (1-x_3/\lambda)},$$ where $\ct_\lambda$ means to take the
constant term in $\lambda$. Then the counting problem is converted
to evaluating the constant term of a special rational function,
which can be done be computer package as in \cite{xiniterate}. For
a rigorous description about how the above works in general
situation, i.e., in a field of iterated Laurent series, the reader
is referred to \cite{xiniterate}.

Using a computer we can easily obtain the generating function of
weak magic squares of order $3$:
\begin{multline*}
  G=  {\frac { \left( 1-tx_{{4}}x_{{7}}x_{{9}}x_{{6}}x_{{2}}x_{{3}}x_{{5}}x_{
{8}}x_{{1}}\right)  \left(
1+tx_{{4}}x_{{7}}x_{{9}}x_{{6}}x_{{2}}x_
{{3}}x_{{5}}x_{{8}}x_{{1}} \right) ^{2}}{ \left(
1-tx_{{1}}x_{{5}}x_{{ 9}}{x_{{4}}}^{2}{x_{{8}}}^{2}{x_{{3}}}^{2}
\right)  \left( 1-tx_{{7}}
x_{{5}}x_{{3}}{x_{{4}}}^{2}{x_{{2}}}^{2}{x_{{9}}}^{2} \right)}}\\
\times \frac{1}{
 \left( 1-tx_{{7}}x_{{5}}x_{{3}}{x_{{1}}}^{2}{x_{{8}}}^{2}{x_{{6}}}^{
2} \right)  \left(
1-tx_{{1}}x_{{5}}x_{{9}}{x_{{7}}}^{2}{x_{{2}}}^{2} {x_{{6}}}^{2}
\right), }
\end{multline*}
where $t$ records $m/3$ since the $m$ is always divisible by $3$,
and the exponents in $x_1$, \dots $x_9$ represents
$a_1,a_2,a_3,b_1,\dots$.

To obtain the generating function for magic squares, we shall take
only terms in $G$ that have different exponents in the $x$'s. To
eliminate those terms with same exponents in $x_1$ and $x_2$, we
subtract by the diagonal $\diag_{x_1,x_2} G$ with respect to $x_1$
and $x_2$, where
$$\diag_{x,y}
\sum_{r\in \NN}\sum_{s\in \NN} b_{r,s} x^ry^s = \sum_{r\in \NN}
b_{r,r} x^r y^r ,$$ and we use the formula for a rational power
series $F(x,y)$:
$$\diag_{x,y} F(x,y)= \ct_{\lambda_1,\lambda_2} \frac{1}{1-xy/(\lambda_1\lambda_2)} F(\lambda_1,\lambda_2).$$
Similarly, we can eliminate those terms with same exponents in
$x_i$ and $x_j$ for all $i$ and $j$.

The generating function of all magic squares of order $3$ is still
complicated. We can add the extra constraints that
$c_3<c_1<a_3<a_1$ to eliminate rotations and reflections. It
suffices to find a way to add the constraint that the exponent of
$x_9$ is smaller than that of $x_7$. The other constraints can be
added iteratively. We omit the details here.

Finally we obtain the generating function of desired magic
squares:
\begin{multline}
 {\frac
{{t}^{4}{x_{{7}}}^{3}{x_{{5}}}^{4}{x_{{3}}}^{5}{x_{{1}}}^{7}{x
_{{8}}}^{8}{x_{{6}}}^{6}x_{{9}}{x_{{4}}}^{2} \left(1+
tx_{{1}}x_{{5}}x_{
{9}}{x_{{4}}}^{2}{x_{{8}}}^{2}{x_{{3}}}^{2}-2{t}^{2}{x_{{5}}}^{2}x_{{9}}
{x_{{4}}}^{2}{x_{{8}}}^{4}{x_{{1}}}^{3}{x_{{3}}}^{3}x_{{7}}{x_{
{6}}}^{2} \right) }{ \left(
1-tx_{{7}}x_{{5}}x_{{3}}{x_{{1}}}^{2}{x_{ {8}}}^{2}{x_{{6}}}^{2}
\right)  \left(1- tx_{{4}}x_{{7}}x_{{9}}x_{{6}}x_
{{2}}x_{{3}}x_{{5}}x_{{8}}x_{{1}}\right)}} \\
\times \frac{1}{ \left( 1-{t}^{2}{x_{{5}}}^{2}x_
{{9}}{x_{{4}}}^{2}{x_{{8}}}^{4}{x_{{1}}}^{3}{x_{{3}}}^{3}x_{{7}
}{x_{{6}}}^{2} \right)  \left(
1-{t}^{3}{x_{{7}}}^{2}{x_{{5}}}^{3}{x_
{{1}}}^{5}{x_{{8}}}^{6}{x_{{3}}}^{4}{x_{{6}}}^{4}x_{{9}}{x_{{4}}}^{2}
 \right) .}\label{e-gf-magicsquare}
\end{multline}
We observe that part of the numerator can be rewritten as
\begin{multline*}
1+tx_{{1}}x_{{5}}x_{
{9}}{x_{{4}}}^{2}{x_{{8}}}^{2}{x_{{3}}}^{2}-2{t}^{2}
{x_{{5}}}^{2}x_{{9}}
{x_{{4}}}^{2}{x_{{8}}}^{4}{x_{{1}}}^{3}{x_{{3}}}^{3}x_{{7}}{x_{
{6}}}^{2} \\
= (1-{t}^{2}{x_{{5}}}^{2}x_{{9}}
{x_{{4}}}^{2}{x_{{8}}}^{4}{x_{{1}}}^{3}{x_{{3}}}^{3}x_{{7}}{x_{
{6}}}^{2})+tx_{{1}}x_{{5}}x_{
{9}}{x_{{4}}}^{2}{x_{{8}}}^{2}{x_{{3}}}^{2}(1-tx_{{7}}x_{{5}}x_{{3}}{x_{{1}}}^{2}{x_{
{8}}}^{2}{x_{{6}}}^{2}),
\end{multline*}
where both terms on the right-hand side will cancel with the
denominator of \eqref{e-gf-magicsquare}. Theorem \ref{t-main} then
follows.

The order $4$ case would be really hard. The difficulty lies in
the fact that there are $880$ pure magic square of order $4$ (up
to rotations and reflections), which  suggests that there will be
at least $880$ simple rational functions. Our current package as
provided in \cite{xiniterate} is not powerful enough to find an
explicit generating function for magic squares of order $4$
analogous to \eqref{e-gf-magicsquare}.

\bibliographystyle{amsplain}

\providecommand{\bysame}{\leavevmode\hbox
to3em{\hrulefill}\thinspace}
\providecommand{\MR}{\relax\ifhmode\unskip\space\fi MR }
\providecommand{\MRhref}[2]{%
  \href{http://www.ams.org/mathscinet-getitem?mr=#1}{#2}
} \providecommand{\href}[2]{#2}

\end{document}